\documentclass[11pt,english]{article}
\usepackage{graphicx} % Requihighlighted for inserting images
\usepackage{amssymb, amsthm, amsmath, fullpage, hyperref}
\usepackage{graphicx, enumerate}
\usepackage{tikz-network}
\usetikzlibrary{positioning, 
                quotes}
\usepackage[utf8]{inputenc}
\usepackage{cleveref}
\usepackage{amsfonts}
\usepackage{amssymb}
\usepackage{amsmath,amsthm,dsfont}
\usepackage[english]{babel}
\usepackage{makeidx}
\usepackage{mathtools}
\usepackage{bbm}
\usepackage{fullpage,amsfonts,amssymb,epsfig,epstopdf,amsmath,titling,url,array,titlesec}
\usepackage{mathtools}
\usepackage{authblk}
\usepackage{tikz}
\usepackage{float}
\usepackage{makeidx}
\usepackage[mathscr]{euscript}
\usepackage{nomencl}
\usepackage{comment}
\theoremstyle{plain}
\newtheorem{thm}{Theorem}[section]
\newtheorem{lem}{Lemma}[section]
\newtheorem{cor}{Corollary}

\theoremstyle{definition}
\newtheorem{defn}{Definition}

\newtheorem{note}[thm]{Note}
\newtheorem{claim}{Claim}

\providecommand{\keywords}[1]
{
  \textbf{\textit{Keywords:}} #1
}

\author[1]{Rajat Adak}
\author[2]{L. Sunil Chandran}
\affil[1]{\texttt{rajatadak@iisc.ac.in}} \affil[2]{\texttt{sunil@iisc.ac.in}}
\affil[1,2]{Department of Computer Science and Automation,
Indian Institute of Science, Bangalore, India}
\date{}

\title{Generalized Zykov's Theorem}

\begin{document}
\maketitle
\begin{abstract}
For a simple graph $G$, let $n$ denote its number of vertices, and let $N(G,K_t)$ denote the number of copies of $K_t$ in $G$. Zykov's theorem (1949) asserts that for any $K_{r+1}$-free graph and $t \ge 2$,
\[
N(G,K_t) \le {r \choose t}\left(\frac{n}{r}\right)^t 
\]
We generalize Zykov’s bound within a vertex-based localization framework

For each vertex $v \in V(G)$, let $c(v)$ denote the order of the largest clique containing $v$. Then
\[
N(G,K_t) \le n^{t-1} \sum_{v \in V(G)} \frac{1}{c(v)^t} {c(v) \choose t}.
\]
Moreover, when $G$ contains a copy of $K_t$, equality holds if and only if $G$ is a regular complete multipartite graph. Note that if we impose the condition that $G$ is $K_{r+1}$-free, then $c(v) \leq r$ for all $v \in V(G)$, and the monotonicity of $s \mapsto \binom{s}{t}/s^t$ gives Zykov's bound.

\end{abstract}
\section{Introduction}
Typical problems in extremal graph theory aim to maximize or minimize certain graph parameters under prescribed structural constraints, and to characterize the graphs that attain these extremal values. Let $m$ and $n$ denote the size and order of the graph, respectively. One of the foundational results in extremal graph theory is the following classical theorem:
\begin{thm}\label{th:Turan}\emph{(Tur\'{a}n \cite{Turan})}
    Let $G$ be a $K_{r+1}$-free graph on $n$ vertices. Then,
    \[
        m \le \frac{n^2(r-1)}{2r},
    \]
    with equality if and only if $G$ is a regular complete $r$-partite graph.
\end{thm}
Let $N(G,K_t)$ denote the number of copies of the clique $K_t$ in the graph $G$. Tur\'{a}n's theorem was later generalized by Zykov~\cite{Zykov1949-qr} in 1949, who established an upper bound on $N(G,K_t)$ for every $t \ge 2$ in any $K_{r+1}$-free graph $G$.

\begin{thm}\label{zykovthm}\emph{(Zykov \cite{Zykov1949-qr})}
    Let $G$ be a $K_{r+1}$-free graph on $n$ vertices. Then for any $t \geq 2$,
    \[N(G,K_t) \leq {r \choose t}\left(\frac{n}{r}\right)^t\]
    Equality holds if and only if $G$ is a regular complete $r$-partite graph.
\end{thm}
\begin{note}
    Note that taking $t = 2$, in \Cref{zykovthm}, we get $N(G,K_t) =m$ and ${r \choose t}\left(\frac{n}{r}\right)^t = \frac{n^2(r-1)}{2r}$. Thus, we can recover Tur\'{a}n's bound from Zykov's bound.
\end{note}
\subsection{Localization}
Recently, Brada\v{c}~\cite{bradac} and Malec and Tompkins~\cite{DBLP:journals/ejc/MalecT23} introduced the notion of \emph{localization} as a framework for strengthening and generalizing classical extremal graph bounds. Traditional extremal results typically depend on global parameters--such as the clique number, the order of the largest clique, as in \Cref{th:Turan} and \ref{zykovthm}.

The localization approach begins by defining a \emph{localization function}. 
For example, the clique number of a graph $G$ is a global parameter, but 
it can be localized by assigning to each edge $e \in E(G)$ a value $w(e)$ equal 
to the size of the largest clique containing $e$. Similarly, the circumference, 
another global parameter, can be localized by assigning to each edge the length 
of the longest cycle containing it (and value $2$ if it is a cut-edge). Many 
other global parameters, such as the length of the longest path, the maximum 
degree, and the girth, can be localized in analogous ways. In this manner, 
localization replaces global constraints with suitable local weight functions 
defined on edges, vertices, or other graph elements.

Brada\v{c}~\cite{bradac} and Malec and Tompkins~\cite{DBLP:journals/ejc/MalecT23} gave an edge-based localization of \Cref{th:Turan};
\begin{thm}\label{thm:GTuran}\emph{(Bradač \cite{bradac} and Malec-Tompkins \cite{DBLP:journals/ejc/MalecT23})} Let $G$ be a graph on $n$ vertices. Then,
\[
    \sum_{e \in E(G)}\frac{w(e)}{w(e) -1} \leq \frac{n^2}{2}
\]
Equality holds if and only if $G$ is a regular complete multipartite graph.
\end{thm}

In \cite{adak2025vertex}, Adak and Chandran initiated the use of the localization framework for generalizing extremal graph problems by assigning weights to vertices. In \cite{adak2025turan}, the same authors developed a vertex-based localization of a \emph{stronger} variant of \Cref{th:Turan} and characterized the extremal graphs.

\begin{defn}
    Let $c: V(G) \rightarrow \mathbb{N}$ be a weight function on the vertices of the graph, such that $c(v)$ is the order of the largest clique containing the vertex $v$, in the graph.
\end{defn}
\begin{thm}\label{th:main}\emph{(Adak-Chandran \cite{adak2025turan})}
    Let $G$ be a graph on $n$ vertices. Then,
\[
    m \leq \left\lfloor\frac{n}{2}\sum_{v\in V(G)}\frac{c(v)-1}{c(v)}\right\rfloor
\]
\end{thm} 
There have been several developments in the localization framework in recent years. In~\cite{zhao2025localized}, Zhao and Zhang applied edge-based localization, for the circumference of the graph, to generalize the Erd\H{o}s--Gallai theorem for cycles~\cite{gallai1959maximal}, while the corresponding path theorem was localized in~\cite{DBLP:journals/ejc/MalecT23}. Kirsch and Nir~\cite{Kirsch_2024} established localized versions of several generalized Tur\'an-type results, including a localization of Zykov's theorem. Arag\~{a}o and Souza~\cite{aragao2024localised} subsequently proved a more general localized graph Maclaurin inequality, which in particular encompasses the localized Zykov result of Kirsch and Nir. 

Adak and Chandran~\cite{adak2025vertexbasedlocalizationgeneralizedturan} developed a vertex-based localization of Luo’s theorems~\cite{luo2018maximum}, which themselves generalize the classical Erd\H{o}s--Gallai bounds. In~\cite{adak2025localizationframeworkgeneralizeextremal}, the same authors provided localized bounds for planar graphs with finite girth, as well as localizations of generalized Tur\'{a}n-type results due to~\cite{wood} and~\cite{CHAKRABORTI2024103955}. 

Furthermore, the localization framework has recently been extended to 
hypergraphs~\cite{DBLP:journals/dm/ZhaoZ25}, and it has also been applied to 
generalize spectral inequalities as in \cite{LIU2026241,liu2025new,kannan2025localizationspectralturantypetheorems}, underscoring its versatility and broad 
applicability within extremal graph theory.

Brada\v{c}~\cite{bradac} adopted the strategy of Nikiforov~\cite{nikiforov2006extension} and Khadzhiivanov~\cite{Khadzhiivanov1977-xg}, which builds on the Motzkin--Straus~\cite{motzkin1965maxima} analytical proof of Tur\'{a}n’s theorem. This approach was later adapted by Arag\~{a}o and Souza~\cite{aragao2024localised}, who obtained a localized version of the graph Maclaurin inequality. Using the same strategy, Liu and Ning developed localized variants of the spectral Tur\'{a}n theorem in~\cite{LIU2026241} and~\cite{liu2025new}.

\section{Our Result}
\begin{thm}\label{ourtheorem} Let $G$ be a graph on $n$ vertices. Then for any $t \geq 2$,
\[N(G,K_t) \leq n^{t-1}\sum_{v\in V(G)}\frac{1}{c(v)^t}{c(v) \choose t}.\]
If $G$ contains a copy of $K_t$, then equality holds if and only if $G$ is a regular complete multipartite graph.
\end{thm}
\begin{note}
If $G$ is $K_t$-free, then both sides of the bound in \Cref{ourtheorem} are zero. Thus the equality characterization is stated only for graphs containing a copy of $K_t$. If $G$ is $K_{r+1}$-free, then $c(v) \leq r$ for all $v \in V(G)$. Since $s \mapsto \binom{s}{t}/s^t$ is increasing for $s \geq t$, \Cref{ourtheorem} gives \Cref{zykovthm}.
\end{note}
\begin{note}The inequality in \Cref{ourtheorem} also follows directly from Theorem~1.2 of Arag\~{a}o and Souza~\cite{aragao2024localised}, by taking $s=1$ and $q=t$, and setting $x_v=a_t(c(v))^{1/t}$. Finally applying  H\"older's inequality yields
 the claimed bound.\end{note} 
 We now give our independent proof.

\subsection{Proof of \Cref{ourtheorem}}

\begin{defn} Let $\Delta^{n-1}$ be the simplex defined as:
\[\Delta^{n-1} = \{x = (x_1,x_2,\dots,x_n) \in \mathbb{R}^n : x_i \geq 0\ \forall\ i \in [n] \text{ and } ||x||_1 =1\}\]
\end{defn}
\begin{defn}
    For a graph $G$ with $n$ vertices and a vector $x \in \Delta^{n-1}$, define the functions $\Phi(G,x), A(G,x), B(G,x)$ such that:
    \[A(G,x) = \sum_{v \in V(G)}\frac{x_v}{c(v)^t}{c(v) \choose t}\]
    \[B(G,x) = \sum_{K_t \subseteq G}\prod_{v \in V(K_t)}x_v\]
    \[\Phi(G,x) = A(G,x) - B(G,x)\]
\end{defn}
\begin{defn}\label{defndelta}
    For $i,j \in V(G)$ and $x \in \Delta^{n-1}$ define;
    \begin{equation}
        \delta_{ij}(x) = \left(\frac{1}{c(i)^t}{c(i) \choose t} - \frac{1}{c(j)^t}{c(j) \choose t}\right) - \left(\sum_{K_{t-1} \subseteq N(i)}\prod_{v \in V(K_{t-1})}x_v -\sum_{K_{t-1} \subseteq N(j)}\prod_{v \in V(K_{t-1})}x_v\right)
    \end{equation}
    Note that for for all $i,j \in V(G)$, $\delta_{ij}(x) = -\delta_{ji}(x)$.
\end{defn}
\begin{defn}
We define a \textit{transfer function} on $\Delta^{n-1}$, which will play a crucial role throughout the paper. 
For $i, j \in V(G)$ and $x \in \Delta^{n-1}$, the transfer function transfers an amount $\epsilon \ge 0$ from the $j$-th coordinate to the $i$-th coordinate of $x$. Formally,
\[
T_{\epsilon}^{ij}(x) = x + \epsilon (e_i - e_j).
\]
Note that if $0\leq \epsilon \leq x_j$, then is it easy to verify that $T_{\epsilon}^{ij}(x) \in \Delta^{n-1}$
\end{defn}
\begin{lem}\label{mainlemma}
    Let $\Phi(G,x) = \text{min}\{\Phi(G,y)\mid y \in \Delta^{n-1}\}$, that is, $x \in \Delta^{n-1}$ be a minimizer of $\Phi(G,\cdot)$. If $i,j \in Supp(x)$ and $i$ is not adjacent to $j$, then $\delta_{ij}(x) = 0$. 
    \begin{proof}
        Let $i, j \in Supp(x)$ such that $i$ is not adjacent to $j$ in $G$. Let $\epsilon \in \mathbb{R}^+$ such that $\epsilon \leq min\{x_i,x_j\}$. Let $y = T_{\epsilon}^{ij}(x)$ and $y'= T_{\epsilon}^{ji}(x)$. 
        
        Let $y= (y_1,y_2,\dots,y_n)$, and $y' = (y'_1,y'_2,\dots,y'_n)$. It is easy to verify that, $y, y' \in \Delta^{n-1}$. Since $x$ is a minimizer, we have: \begin{equation}\label{y>x}
            \Phi(G,y) - \Phi(G,x) \geq 0
        \end{equation}
        Since $\Phi(G,\cdot) = A(G,\cdot) - B(G,\cdot)$, it is easy to see that:
        \begin{equation}\label{eq1}
                        \Phi(G,y) - \Phi(G,x) = (A(G,y) - A(G,x)) - (B(G,y) - B(G,x))
        \end{equation}
        Since $y_i = x_i +\epsilon$, $y_j = x_j - \epsilon$, and $y_k =x_k$ for all $k \notin \{i,j\}$ we have;
\begin{align}
    A(G,y) - A(G,x) &= \frac{y_i-x_i}{c(i)^t}{c(i) \choose t} + \frac{y_j - x_j}{c(j)^t}{c(j) \choose t}\nonumber\\
    &= \epsilon \left(\frac{1}{c(i)^t}{c(i) \choose t} - \frac{1}{c(j)^t}{c(j) \choose t}\right)\label{eq2}
\end{align}

Because vertices $i$ and $j$ are non-adjacent, no copy of $K_t$ in $G$ can contain both simultaneously. Consequently, every $K_t$ must contain either exactly one of these vertices or neither. Observe that those copies disjoint from $\{i, j\}$ cancel out and do not affect the difference $B(G,y_{\epsilon}) - B(G,x)$. 
\begin{align}
    B(G,y) - B(G,x) &= \sum_{\substack{K_t \subseteq G\\ i \in V(K_t)}}\prod_{v\in V(K_t)}y_v - \sum_{\substack{K_t \subseteq G\\ i \in V(K_t)}}\prod_{v\in V(K_t)}x_v  + \sum_{\substack{K_t \subseteq G\\ j \in V(K_t)}}\prod_{v\in V(K_t)}y_v - \sum_{\substack{K_t \subseteq G\\ j \in V(K_t)}}\prod_{v\in V(K_t)}x_v\nonumber\\
    &=(y_i - x_i)\sum_{K_{t-1} \subseteq N(i)}\prod_{v \in V(K_{t-1})}x_v + (y_j - x_j)\sum_{K_{t-1} \subseteq N(j)}\prod_{v \in V(K_{t-1})}x_v\nonumber\\
    &=\epsilon\left(\sum_{K_{t-1} \subseteq N(i)}\prod_{v \in V(K_{t-1})}x_v -\sum_{K_{t-1} \subseteq N(j)}\prod_{v \in V(K_{t-1})}x_v\right)\label{eq3}
\end{align}
From \cref{defndelta} and using \cref{eq1,eq2,eq3} we get that $\Phi(G,y) - \Phi(G,x) = \epsilon\delta_{ij}(x)$. Thus from \cref{y>x} we get $\epsilon\delta_{ij}(x) \geq 0$. 

Note that, since $x$ is a minimizer, and $y' \in \Delta^{n-1}$, we also have $\Phi(G,y') -\Phi(G,x)\geq 0$. Thus using similar analysis as above we get that $\Phi(G,y') -\Phi(G,x) = \epsilon\delta_{ji} \geq 0$.

From \cref{defndelta} we have $\delta_{ij} = -\delta_{ji}$. Therefore, $\epsilon\delta_{ji} \geq 0 \implies \epsilon\delta_{ij} \leq 0$. But we already have $\epsilon\delta_{ij} \geq 0$. Thus, $\epsilon\delta_{ij}=0$, and since $\epsilon > 0$, we have $\delta_{ij} = 0$.
    \end{proof}
\end{lem}
\begin{cor}\label{cor1}
Let $x \in \Delta^{n-1}$ be a minimizer of $\Phi(G,\cdot)$. If $i,j \in \operatorname{Supp}(x)$ and $i$ is not adjacent to $j$, then both $y = T_{x_j}^{ij}(x)$ and $y' = T_{x_i}^{ji}(x)$ are also minimizers of  $\Phi(G,\cdot)$.

    \begin{proof}
        Following similar steps as in \Cref{mainlemma} we get that, $\Phi(G,y) - \Phi(G,x) = x_j\delta_{ij}$ and $\Phi(G,y') - \Phi(G,x) = x_i\delta_{ji}$. But $\delta_{ij} = \delta_{ji} = 0$, since $i$ and $j$ are not adjacent. Therefore, $\Phi(G,y) - \Phi(G,x) = \Phi(G,y') - \Phi(G,x) = 0$. Thus, $y$ and $y'$ are minimizers of $\Phi(G,\cdot)$.
    \end{proof}
\end{cor}
\begin{cor}\label{cor2}
    Let $x \in \Delta^{n-1}$ be a minimizer of $\Phi(G,\cdot)$. If $S \subseteq \operatorname{Supp}(x)$ is an independent set of $G$, 
then transferring the total amount from all coordinates indexed by $S$ 
to the coordinate of a single vertex in $S$ produces another minimizer. Therefore, if $i \in S$, we can construct another minimizer $z = x + \sum_{t \in S \setminus \{i\}} x_t(e_i - e_t)$ by shifting the total amount in $S$ to the vertex $i$.

    \begin{proof}
This follows by repeatedly applying \cref{cor1} to $(i,t)$ with 
$t \in S \setminus \{i\}$, transferring the remaining weight from $t$ to $i$. 
In this way, all the weight on the coordinates indexed by $S$ is successively 
merged onto the $i^{\text{th}}$ coordinate, while preserving minimality throughout.

    \end{proof}
\end{cor}
\begin{cor}\label{cor3}
    If $x \in \Delta^{n-1}$ is a minimizer of $\Phi(G,\cdot)$ with minimal support, then $G[Supp(x)]$ is a clique.
    \begin{proof}
        Suppose not, then there exists $i,j \in Supp(x)$ such that $i$ is not adjacent to $j$. Thus from \cref{cor1}, $y = T_{x_j}^{ij}(x)$ is a minimizer. Note that $y_i = x_i + x_j$, $y_j=0$ and $y_k = x_k$ for all $ k \notin \{i,j\}$. Thus $Supp(y) \subsetneq Supp(x)$. This contradicts the minimality of $Supp(x)$. Thus $Supp(x)$ induces a clique in $G$.
    \end{proof}
\end{cor}

\subsubsection{Proof of Inequality}
Now we are all set to prove the inequality of \Cref{ourtheorem}.
\begin{proof}
    Let $x$ be a minimizer of $\Phi(G,\cdot)$ with minimal support. From \cref{cor3} we obtain that $G[Supp(x)]$ is a clique. Let $|Supp(x)| = k$. Clearly, $c(v) \geq k$ for all $v \in Supp(x)$. It is easy to see that $f(x) = \frac{1}{x^t}{x \choose t}$ is an increasing function for $x \geq 1$. Thus, for all $v\in Supp(x)$ we have:
    \begin{equation}\label{eqn4}
        \frac{1}{c(v)^t}{c(v) \choose t} \geq \frac{1}{k^t}{k \choose t}
    \end{equation}
 Since for all $i \notin Supp(x)$, $x_i = 0$ we get that:
    \begin{equation}\label{eqn5}
        B(G,x)=\sum_{K_t \subseteq V(G)}\prod_{v \in V(K_t)}x_v = \sum_{K_t \subseteq Supp(x)}\prod_{v \in V(K_t)}x_v
    \end{equation}
Now since $\sum_{v \in Supp(x)}x_v = 1$, using Maclaurin's inequality we get that:
\begin{align}
    \frac{\sum_{v \in Supp(x)}x_v}{k} &\geq \left(\frac{\sum_{K_t \subseteq Supp(x)}\prod_{v \in V(K_t)}x_v}{{k\choose t}}\right)^{\frac{1}{t}}\nonumber\\
    \implies \frac{1}{k^t}{k\choose t}  &\geq \sum_{K_t \subseteq Supp(x)}\prod_{v \in V(K_t)}x_v\label{eqn6}
\end{align}
From \cref{eqn5,eqn6} we get that:
\begin{equation}\label{eqn7}
    B(G,x) \leq \frac{1}{k^t}{k\choose t}
\end{equation}
Thus combining \cref{eqn4,eqn7} we get that:
\begin{align}
    \Phi(G,x) &= \sum_{v \in V(G)}\frac{x_v}{c(v)^t}{c(v) \choose t} - \sum_{K_t \subseteq G}\prod_{v \in V(K_t)}x_v\nonumber\\
    &\geq \sum_{v \in Supp(x)}\frac{x_v}{k^t}{k \choose t} - \frac{1}{k^t}{k\choose t}\nonumber\\
    &=\frac{1}{k^t}{k\choose t}-\frac{1}{k^t}{k\choose t} = 0 \label{eqn8}
\end{align}
Since $x$ is a minimizer, we get that $\Phi(G,y) \geq 0$ for all $y \in \Delta^{n-1}$

Let $y = (\frac{1}{n}, \frac{1}{n}, \dots, \frac{1}{n})$. Clearly $y \in \Delta^{n-1}$. Therefore,
\begin{align}
    \Phi(G,y) &= \frac{1}{n}\sum_{v \in V(G)}\frac{1}{c(v)^t}{c(v) \choose t} - \sum_{K_t \subseteq G}\prod_{v \in V(K_t)}\frac{1}{n}\nonumber\\
    &= \frac{1}{n}\sum_{v \in V(G)}\frac{1}{c(v)^t}{c(v) \choose t} - \frac{1}{n^t}N(G,K_t)\geq 0\label{eqn9}\\
    &\implies N(G,K_t) \leq n^{t-1}\sum_{v \in V(G)}\frac{1}{c(v)^t}{c(v) \choose t}\nonumber
\end{align}
Thus, we get the inequality as in \Cref{ourtheorem}.
\end{proof}
\subsubsection{Characterizing Extremal Graphs}
Since $t \geq 2$, it is easy to check that isolated vertices do not contribute to either side of the bound in \Cref{ourtheorem}. Thus we can assume $G$ does not contain any isolated vertices.
\vspace{2mm}
\newline\textbf{$\bullet$ If Part}
\newline First, we will show that if $G$ is a regular complete multipartite graph, then we get equality in the bound of \Cref{ourtheorem}.

Let $G$ be a regular complete $r$-partite graph. Clearly $c(v) = r$ for all $v \in V(G)$, therefore $G$ is $K_{r+1}$-free. The right-hand side of the bound in \Cref{ourtheorem} coincides with the upper bound of Zykov's theorem. Since regular complete $r$-partite graphs are extremal for \Cref{zykovthm}, we get equality in \Cref{ourtheorem}. 

\vspace{2mm}
\noindent\textbf{$\bullet$ Only If Part}
\newline For the other way, we will assume that we have equality in the bound of \Cref{ourtheorem}. Thus, we must have equality in \cref{eqn8,eqn9}. Thus, we get the following \textbf{equality conditions}:
\begin{enumerate}
    \item\label{ec3} From equality in \cref{eqn9} we have $\Phi(G,y) = 0$, where $y = \{\frac{1}{n},\frac{1}{n},\dots,\frac{1}{n}\}$. Thus $y$ is a minimizer of $\Phi(G,\cdot)$.

    From the above condition we know that minimum of $\Phi(G,\cdot)$ is $0$, then for any minimizer $x \in \Delta^{n-1}$, $\Phi(G,\cdot) = 0$. Thus if $x$ is a minimal support minimizer, from \cref{eqn8} we get \[\sum_{v \in V(G)}\frac{x_v}{c(v)^t}{c(v) \choose t} - \sum_{K_t \subseteq G}\prod_{v \in V(K_t)}x_v = 0\] Thus we must have equality in \cref{eqn4,eqn6}, that is, \[\frac{1}{c(v)^t}{c(v) \choose t} = \frac{1}{k^t}{k \choose t}\]  \[\frac{1}{k^t}{k\choose t}  = \sum_{K_t \subseteq Supp(x)}\prod_{v \in V(K_t)}x_v\]
    
    \item\label{ec1} Equality holds in Maclaurin inequality if and only if all the entries are equal. Thus from the equality of \cref{eqn6} we get that, if $x$ is a minimal support minimizer of $\Phi(G,\cdot)$ then, for all $v \in Supp(x)$, $x_v = \frac{1}{|Supp(x)|}$.
    \item\label{ec2} From equality in \cref{eqn4}, we have that, if $x$ is a minimal support minimizer, then for all $v \in Supp(x)$, $c(v) = |Supp(x)|$ 
\end{enumerate}
\begin{lem}\label{auxlemma}
    If $z \in \Delta^{n-1}$ is a minimizer of $\Phi(G,\cdot)$ and $G[Supp(z)]$ is a clique then $z$ is a minimal support minimizer.
    \begin{proof}
        Suppose $z$ is not a minimal support minimizer, then there exists a minimizer $z' \in \Delta^{n-1}$ such that $Supp(z')\subsetneq Supp(z)$ and $z'$ has minimal support. From the equality condition~\ref{ec2}, we get that $c(v) = |Supp(z')| <|Supp(z)|$ for all $v \in Supp(z')$. But by assumption the vertices in $Supp(z')$ are already part of a $|Supp(z)|$ order clique. Thus, we get a contradiction.
    \end{proof}
\end{lem}
\begin{lem}\label{mainlemma2}
    If $z \in \Delta^{n-1}$ is a minimizer of $\Phi(G,\cdot)$, then $G[Supp(z)]$ is a complete $\omega_z$-partite graph, where $\omega_z$ is the clique number of $G[Supp(z)]$. Furthermore, if $V_1,V_2,\dots, V_{\omega_z}$ are the parts of $G[Supp(z)]$, then for all $k \in [\omega_z]$, $\sum_{v \in V_k}z_v = \frac{1}{\omega_z}$.
    \begin{proof}
Suppose there exists a minimizer of $\Phi(G,\cdot)$ in $\Delta^{n-1}$ for which the above
statement fails. Let $z$ be such a minimizer with minimal support; that is, for every 
$z' \in \Delta^{n-1}$ with $\operatorname{Supp}(z') \subsetneq \operatorname{Supp}(z)$, 
the lemma holds for $z'$.

If $G[Supp(z)]$ is a clique then from \Cref{auxlemma} we get that $z$ is a minimal support minimizer. Note that, $\omega_z = |Supp(z)|$ and $G[Supp(z)]$ is complete $\omega_z$-partite. Also from the equality condition~\ref{ec1} we have $z_v = \frac{1}{|Supp(z)|} = \frac{1}{\omega_z}$ for all $v \in Supp(z)$. Thus we get a contradiction.

Thus assume $G[Supp(z)]$ is not a clique. Therefore, there exists $i,j \in Supp(z)$ such that $i$ and $j$ are not adjacent. From \cref{cor1}, we have $\Phi(G,y) = \Phi(G,z)$, where $y = T_{z_j}^{ij}(z)$. Note that $Supp(y) = Supp(z)\setminus \{j\}$, therefore $Supp(y) \subsetneq Supp(z)$. Since, $y$ is a minimizer of $\Phi(G,\cdot)$, from our assumption we get that $y$ satisfies the lemma and therefore $G[Supp(y)]$ is a complete $\omega_{y}$-partite graph. Let $P_1,P_2,\dots,P_{\omega_y}$ be the parts in $G[Supp(y)]$.
\begin{claim}\label{claim1}
            There exists $k \in [\omega_y]$ such that $N(j)\cap P_k = \emptyset$.
            \begin{proof} Assume on the contrary, there exists $v_k \in N(j)\cap P_k$ for all $k \in [\omega_y]$. Note that we can define $y^* = (y^*_1,y^*_2,\dots,y^*_n) \in \Delta^{n-1}$ such that $y^*_{v_k} = \sum_{v \in P_k}y_v$ for all $k\in [\omega_y]$ and $y^*_l = 0$ for all $l \notin \{v_1,v_2,\dots,v_{\omega_y}\}$ (that is, for $k \in [\omega_y]$, we transfer all the amount from the coordinates indexed by $P_k$ to the coordinate corresponding $v_k$). Therefore, $Supp(y^*) = \{v_1,v_2,\dots v_{\omega_y}\}.$ Since each $P_k$ is an independent set, from \cref{cor2}, we get that $\Phi(G,y^*) = \Phi(G,y) $. Thus, $y^*$ is a minimizer. Clearly, $Supp(y^*)$ induces an $\omega_y$-order clique. Therefore, from \Cref{auxlemma} we get that $y^*$ is a minimal support minimizer. From equality condition~\ref{ec2} we have $c(v_k) = |Supp(y^*)| = \omega_y$. But from the assumption we get that $Supp(y^*)\cup \{j\}$ induces an $\omega_y+1$-order clique. Since $v_k$ is contained in this clique, we have $c(v_k) \geq \omega_{y}+1$, which results a contradiction.
            \end{proof}
        \end{claim}
        Without loss of generality assume $i \in P_1$.
        \begin{claim}\label{claim2}
            For $k>1$, $N(j) \cap P_{k} \neq \emptyset$.
            \begin{proof}
                Suppose there exists $k >1$, such that $P_k \cap N(j) = \emptyset$.
                \vspace{1mm}
                \newline First suppose $|P_1|>1$. Let $y' = T_{z_i}^{ji}(z)$. Clearly from \cref{cor1}, $\Phi(G,y') = \Phi(G,z)$, that is, $y'$ is also a minimizer. Moreover, $Supp(y')\subsetneq Supp(z)$. Thus, from our assumption, the lemma holds for $y'$. Therefore, $G[Supp(y')|$ is a complete $\omega_{y'}$-partite graph. Clearly $P'_1 =P_1\setminus\{i\}, P_2,\dots, P'_k =P_k\cup\{j\},\dots, P_{\omega_y}$ are maximal independent sets in $G[Supp(y')]$. It is easy to see that $\{P'_1,P_2,\dots,P'_k,\dots,P_{\omega_y}\}$ form the parts in $G[Supp(y')]$, therefore $\omega_{y'} = \omega_y$. Note that $y'_v = y_v$ for all $v \notin \{i,j\}$, $y'_i = y_j = 0$, and $y'_j = y_i = z_i+z_j$. As $y$ satisfies the lemma, $\sum_{v \in P_1}y_v = \sum_{v \in P_k}y_v =\frac{1}{\omega_y}$. But then we get $\sum_{v \in P'_1}y'_v < \sum_{v \in P'_k}y'_v$, which results in a contradiction, since the lemma holds for $y'$ and therefore,  $\sum_{v \in P'_1}y'_v = \sum_{v \in P'_k}y'_v$.
                \vspace{1mm}
                \newline Now suppose $|P_1| =1$. Clearly $\{P_2,\dots,P'_k,\dots,P_{\omega_y}\}$ are the parts in $G[Supp(y')]$, where $P'_k = P_k\cup\{j\}$.  Suppose that $G[\operatorname{Supp}(y')]$ is $1$-partite, that is, $P_k'$ is the only part.  
                Using \cref{cor2}, by transferring the total amount from the coordinates indexed by $P_k'$ to the coordinate corresponding to a single vertex $v \in P_k'$, we obtain another minimizer whose support consists solely of $v$.  
                Then, by equality condition~\ref{ec2}, we have $c(v)=1$, therefore, $v$ is an isolated vertex. However, this contradicts our assumption that $G$ has no isolated vertices.
                
                Thus, $P'_k$ is not the only part in $G[Supp(y')]$. Since the lemma holds for $y$, we have $\sum_{v \in P_t}y_v = \sum_{v \in P_k}y_v$, where $t \notin \{1,k\}$. Note that, $\sum_{v \in P_t}y_v = \sum_{v \in P_t}y'_v$ and $\sum_{v \in P_k}y_v < \sum_{v \in P'_k}y'_v$. Thus we get, $\sum_{v \in P_t}y'_v < \sum_{v \in P'_k}y'_v$, which results in a contradiction, since $y'$ must also satisfy the lemma.
            \end{proof}
        \end{claim}
        From \cref{claim1}, there exists $k\in[\omega_y]$ such that $N(j)\cap P_k = \emptyset$, and from \cref{claim2}, for $k>1$, $N(j)\cap P_k \neq \emptyset$. Thus $N(j) \cap P_1 = \emptyset$.
        \begin{claim}\label{claim3}
             For $k>1$, $N(j)\cap P_k = P_k$.
            \begin{proof}
                               
                Suppose for some $k>1$, there exists $v_k \in P_k$ such that $v_k$ is not adjacent to $j$. But from \cref{claim2}, $P_k \cap N(j) \neq \emptyset$. Thus there exists $u_k \in P_k$ which is adjacent to $j$. In $G[Supp(y')]$, where $y' = T_{x_i}^{ji}(z)$, $v_k$ and $u_k$ are in the same part but $N(v_k) \neq N(u_k)$. This contradicts the fact that $G[Supp(y')]$ is a complete multipartite graph.
            \end{proof}
        \end{claim}
        Thus from \cref{claim3} we get that $j$ is adjacent to all the vertices in $\bigcup_{i=2}^{\omega_y}P_i$. Also we know that $N(j)\cap P_1 = \emptyset$. Therefore, $G[Supp(z)]$ is a complete multipartite graph, with parts $V_1 = P_1\cup\{j\}, V_2 = P_2,\dots, V_{\omega_y}=P_{\omega_y}$. Note that $\omega_z = \omega_y$. 
        
        Since $y$ satisfies the lemma, we know that, for all $k \in [\omega_z]$, $\sum_{v \in P_k}y_v = \frac{1}{\omega_y}$. Note that $i,j \in V_1$ and we have, $y_v = z_v$ for all $v \notin \{i,j\}$. Therefore, for all $k>1$ we get,
        \[\sum_{v \in V_k}z_v = \sum_{v \in P_k}y_v = \frac{1}{\omega_y} = \frac{1}{\omega_z}\]
        Since $\sum_{v \in V(G)}z_v = \sum_{v \in Supp(z)}z_v = 1$, we get 
        \[\sum_{v \in V_1} z_v= 1 - \sum_{k =2}^{\omega_y}\left(\sum_{v \in V_k}z_v\right) = 1 - (\omega_y -1)\frac{1}{\omega_y} = \frac{1}{\omega_z}\] 
        Thus, for all $k \in [\omega_z]$, \[\sum_{v \in V_k}z_v  = \frac{1}{\omega_z}\]
        Thus, $z$ satisfies \Cref{mainlemma2}, which contradicts our initial assumption.
    \end{proof}
\end{lem}
From the equality condition~\ref{ec3}, we obtain that 
$ y = \left\{\frac{1}{n}, \frac{1}{n}, \dots, \frac{1}{n}\right\} $ 
is a minimizer of $ \Phi(G, \cdot) $. Here $Supp(y) = V(G)$. Hence, by \Cref{mainlemma2}, the graph $G$ must be complete multipartite. Let $V_p$ and $V_q$ be any two parts of $G$. Using \Cref{mainlemma2} again, we have
\[
\sum_{v \in V_p} \frac{1}{n} = \sum_{v \in V_q} \frac{1}{n}
\quad \Longrightarrow \quad
\frac{|V_p|}{n} = \frac{|V_q|}{n} \implies |V_p| = |V_q|
\]
Since $V_p$ and $V_q$ were chosen arbitrarily, it follows that all parts of $G$ have the same size.
Thus, $G$ is a regular complete multipartite graph.

\bibliographystyle{plain}
\bibliography{references}
\end{document}